\documentclass[11pt]{article}
\usepackage[dvips]{graphics}
\usepackage{multirow}
\usepackage{array}
\usepackage{epsfig,rotating}
\usepackage{amssymb,amsmath}
\usepackage{latexsym}
\usepackage[usenames]{color}

\numberwithin{equation}{section}

\newcommand{\be}{\begin{equation}}
\newcommand{\ee}{\end{equation}}
\newcommand{\beaa}{\begin{eqnarray*}}
\newcommand{\eeaa}{\end{eqnarray*}}
\newcommand{\bea}{\begin{eqnarray}}
\newcommand{\eea}{\end{eqnarray}}

\newcommand{\bei}{\begin{itemize}}
\newcommand{\eei}{\end{itemize}}

\newtheorem{theorem}{\noindent T{\footnotesize HEOREM}}

\newtheorem{lemma}{\noindent L{\footnotesize EMMA}}[section]
\newtheorem{coro}{\noindent C{\footnotesize OROLLARY}}[section]

\newcommand{\bz}{\mathbf{z}}
\begin{document}

\title{\Large Spectral Radii of Products of Random Rectangular Matrices}

\author{Yongcheng Qi$^1$~~~and~~~ Mengzi Xie$^2$\\
University of Minnesota Duluth
}

\date{}
\maketitle

\footnotetext[1]{Department of Mathematics and Statistics, University of Minnesota Duluth,
MN 55812, USA. {Email: yqi@d.umn.edu (corresponding author)}.
}
\footnotetext[2]{Department of Mathematics and Statistics, University of Minnesota Duluth,
MN 55812, USA. {Email: xie00035@d.umn.edu}.
}




\begin{abstract}

\noindent We consider $m$ independent random rectangular matrices whose entries are independent and identically distributed standard complex Gaussian random variables. Assume the product of the $m$ rectangular matrices is an $n$ by $n$ square matrix. The maximum absolute values of the $n$ eigenvalues of the product matrix is called spectral radius. In this paper, we study the limiting spectral radii of the product when $m$ changes with $n$ and can even diverge.  We give a complete description for the limiting distribution of the spectral radius.  Our results reduce to those in Jiang and Qi \cite{JiangQi2017} when the rectangular matrices are square ones.
\end{abstract}

\noindent \textbf{Keywords:\/} spectral radius, eigenvalue, random rectangular matrix,  non-Hermitian random matrix.

\noindent\textbf{AMS 2010 Subject Classification: \/} 15B52, 60F99, 60G70, 62H10.\\

\newpage
\section{Introduction}\label{intro}

Since Wishart's \cite{Wishart} work on large covariance matrices in multivariate analysis,  the study of random matrices has drawn much attention from mathematics
and physics communities and has found applications in areas such as  heavy-nuclei
(Wigner \cite{Wigner}), condensed matter physics (Beenakker \cite{Been1997}), number theory (Mezzadri and Snaith \cite{MS2005}), wireless communications (Couillet and Debbah \cite{CD}), and high dimensional statistics (Johnstone \cite{John2001, John2008}, and Jiang \cite{Jiang09}).  Bouchaud and Potters~\cite{bouchaud2009financial} provide a survey on applications in finance. The interested reader can find more references in the Oxford Handbook of Random Matrix Theory by Akemann, Baik and Francesco \cite{akemann2011oxford}.

Random matrix theory studies the eigenvalues of random matrices, including the properties of the spectral radii and the empirical spectral distributions of the eigenvalues.  Tracy and Widom \cite{Tracy94, Tracy96} show that the largest eigenvalues of the three Hermitian matrices (Gaussian orthogonal ensemble, Gaussian unitary ensemble and Gaussian symplectic ensemble) converge in distribution to some limits which are  now  known as Tracy-Widom laws. Subsequently,
the Tracy-Widom laws have found more applications, see, e.g., Baik et al. \cite{Baik1999}, Tracy and Widom~\cite{TW02}, Johansson~\cite{Johansson07},  Johnstone~\cite{John2001, John2008} and Jiang \cite{Jiang09}.

The study of non-Hermitian matrices has also attracted attention in the literature. Theoretical results in this direction can be applied to quantum chromodynamics, choaotic quantum systems and growth processes, dissipative quantum maps and  fractional quantum Hall effect.  More applications can be found in Akemann et al.~\cite{akemann2011oxford} and Haake~\cite{haake1991quantum}.
 In the stimulating work by Rider~\cite{rider2003limit, rider2004order} and Rider and Sinclair~\cite{rider2014extremal},
 the spectral radii of the real, complex and symplectic Ginibre ensembles are investigated. It is shown that the spectral radius of the complex Ginibre ensemble   converges to the Gumbel distribution. Jiang and Qi~\cite{JiangQi2017} study the largest radii of three rotation-invariant and non-Hermitian random matrices: the spherical ensemble, the truncation of circular unitary ensemble and the product ensemble, and Jiang and Qi~\cite{JiangQi2018} investigate the limiting empirical spectral distributions for two types of product ensembles. More related work can be also found in Gui and Qi~\cite{GuiQi2018},  Chang and Qi~\cite{ChangQi2017},
Chang, Li and Qi~\cite{ChangLiQi2018}, and Zeng~\cite{Zeng2016, Zeng2017}. The study of the lower and upper tail probabilities of the largest radii is also
of interest, see, e.g.,  Lacroix-A-Chez-Toine et al.~\cite{La2018} and references therein.

Products of random matrices are particularly of interest in recent research. Ipsen~\cite{Ipsen2015} provides several applications, include
wireless telecommunication, disordered spin chain, the stability of large complex system, quantum transport in disordered wires, symplectic maps and Hamiltonian mechanics, quantum chromo-dynamics at non-zero chemical potential. Here we will do a very brief survey for recent developments on the limiting spectral radii and empirical spectral distributions for product ensembles. Two recent papers by Jiang and Qi~\cite{JiangQi2017, JiangQi2018} consider the spectral radii and empirical spectral distribution for the product of $m$ independent $n$ by $n$ Ginibre ensembles, where $m$ can change with $n$ and obtain the limiting distribution functions for the spectral radii and limiting empirical spectral distributions. For earlier works on empirical spectral distribution for
the product ensembles for fixed $m$, see, e.g.,  G\"{o}tze and Tikhomirov \cite{Goetz} , Bordenave~\cite{Bor}, O'Rourke and Soshnikov~\cite{Rourke}, O'Rourke et al.~\cite{Rourke14}, Burda et al.~\cite{BJW},  Burda~\cite{Burda}, and Bai~\cite{Bai}.  Jiang and Qi~\cite{JiangQi2018} also investigate the limiting empirical spectral distribution for the product of $m$ independent truncated Haar unitary matrices when $m$ changes with the dimension of the product matrices.
For the products of $m$ independent spherical ensembles, Chang, Li and Qi~\cite{ChangLiQi2018} study the limiting spectral radius when $m$ can change with the dimension of the product matrices,  Zeng~\cite{Zeng2017} and Chang and Qi~\cite{ChangQi2017} investigate the empirical spectral distribution for the products.

In this paper, we consider the product of $m$ random rectangular matrices with independent and identically distributed (i.i.d.) complex Gaussian entries and investigate the limiting distributions for
the spectral radii. When $m$ is a fixed integer, Zeng~\cite{Zeng2017} obtains the limiting empirical spectral distribution. When these rectangular matrices are actually squared ones, the product matrix is reduced to the product of Ginibre ensembles, which has been studied in Jiang and Qi~\cite{JiangQi2017}.
The products of rectangular matrices have found applications in wireless telecommunication and econophysics (Akemann et al.~\cite{akemann2013products}, Muller~\cite{muller2002asymptotic}, Tulino and Verd~\cite{tulino2004random}),  transport in disordered and chaotic dynamical system(Crisanti et al.~\cite{crisanti1994products}, Ipsen and Kieburg~\cite{ipsen2014weak}).  In particular,  for $m=2$, the product can be regarded as the asymmetric correlation matrices (Vinayak~\cite{vinayak2013spectral}, Vinayak and Benet~\cite{benet2014spectral}) and  has been widely used in finance (Bouchaud et al.~\cite{bouchaud2007large}, Bouchaud and Potters~\cite{bouchaud2009financial}, Livan and Rebecchi~\cite{livan2012asymmetric}).

The rest of the paper is organized as follows. In Section~\ref{main}, we introduce the main results of the paper. In Section~\ref{proofs}, we present some preliminary lemmas and give the proofs for the main results.


\section{Main Results}\label{main}

For integer $m\ge 1$, assume  $\{n_r, ~1\le r\le m+1\}$ are positive integers such that $n_{1}=n_{m+1}=\min\{n_1, \cdots, n_{m+1}\}$. Write $n=n_1=n_{m+1}$ for convenience.  For each $r\in \{1,\cdots, m\}$, $A_r$ is an $n_r\times n_{r+1}$ random rectangular matrix given by
$$
A_{r}={\left(
\begin{array}{cccc}
g_{11}^{(r)}&g_{12}^{(r)}&\cdots & g_{1n_{r+1}}^{(r)}\\
g_{21}^{(r)}&g_{22}^{(r)}&\cdots&g_{2n_{r+1}}^{(r)} \\
\vdots & \vdots & \ddots & \vdots \\
g_{n_{r}1}^{(r)}&g_{n_{r}2}^{(r)}& \cdots &g_{n_{r}n_{r+1}}^{(r)}
\end{array}
\right)},
$$
where $g_{ij}^{(r)}$, $1\leq i\leq n_{r}$, $1\leq j\leq n_{r+1}$ are i.i.d. standard complex normal random variables with $\mathbb{E}g_{ij}^{(r)}=0$, $\mathbb{E}|g_{ij}^{(r)}|^{2}=1$ for $1\le i\le n_r$, $1\le j\le n_{r+1}$, $r=1,\cdots, m$.

Define $A_n^{(m)}$ as the product of the $m$ rectangular matrices $A_{r}$'s, that is,  $A_{n}^{(m)}=A_{1}\cdots A_{m}$.
Let $\bz_{1},\cdots,\bz_{n}$ be the eigenvalues of $A_{n}^{(m)}$.
Set $l_{r}=n_{r}-n$, $r=1,\cdots, m$.
The joint density function for $\bz_{1},\cdots,\bz_{n}$, given in Theorem 2 of Adhikari~\cite{adhi2016}, is as follows
\begin{equation} \label{eigen}
 p(z_{1},\cdots,z_{n})=C\prod\limits_{1\leq j<k\leq n}{\left|z_{j}-z_{k}\right|^{2}}\prod\limits_{j=1}^{n}{w_{m}^{(l_{1},\cdots,l_{m})}}{(\left|z_{j}\right|)}
\end{equation}
with respect to the Lebesgue measure on $\mathbb{C}^{n}$, where $C$ is a normalizing constant, and function $w_{m}^{(l_{1},\cdots,l_{m})}(z)$ can be obtained
recursively by
\[
w_{k}^{(l_{1},\cdots,l_k)}(z)=2\pi\int_{0}^{\infty}w_{k-1}^{(l_{1},\cdots,l_{k-1})}(\frac{z}{s})w_{1}^{(l_k)}(s)\frac{ds}{s}, ~~k\ge 2
\]
with initial
$w_{1}^{(l)}(z)=\exp(-\left|z\right|^{2})\left|z\right|^{2l}$ for any $z$ in the complex plane (see, Zeng~\cite{Zeng2017}).


The spectral radius of $A_{n}^{(m)}$ is defined as the maximal absolute value of the $n$ eigenvalues $\bz_1, \cdots, \bz_n$, i.e.  $\max\limits_{1\le j\le n}|\bz_j|$. In this paper we aim at the limiting distribution of $\max\limits_{1\leq j\leq n}|\bz_{j}|$.  We allow that $m$ changes with $n$. From now on we will write $m$ as $m_n$.


We need to define some notation before we introduce the main results.

Define $\Phi(x)=\frac{1}{\sqrt{2\pi}}\int^x_{-\infty}e^{-t^2/2}dt$ as the standard normal cumulative distribution function (cdf) and $\Lambda(x)=\exp(-e^{-x})$ as the Gumbel distribution function.
For $\alpha\in(0,\infty)$, set
$$\Phi_{\alpha}(x)=\prod\limits_{j=0}^{\infty}\Phi(x+j\alpha^{1/2}),$$  $\Phi_{0}(x)=\Lambda(x)=\exp(-e^{-x})$, and $\Phi_{\infty}(x)=\Phi(x)$. The digamma function $\psi$ is defined by
\begin{equation}\label{psi}
\psi(z)=\frac{d}{dz}\ln\Gamma(z)=\frac{\Gamma'(z)}{\Gamma(z)},
\end{equation}
where $\Gamma(z)$ is the Gamma function.
For large $y$, define
\begin{equation}\label{ab}
a(y)=(\ln y)^{1/2}-(\ln y)^{-1/2}\ln(\sqrt{2\pi}\ln y)~~\mbox{  and }~~ b(y)=(\ln y)^{-1/2}.
\end{equation}

Now we define
\[
\Delta_n=\sum\limits_{r=1}^{m_n}\frac{1}{n_{r}}.
\]
The limiting spectral radius depends on the limit of $\Delta_n$.

We first give a general result on the limiting distribution for the logarithmic spectral radii.


\begin{theorem} \label{t1}
Assume that $\bz_{1},\cdots,\bz_{n}$ are the eigenvalues of $A_{n}^{(m_{n})}$, and
\begin{equation}\label{alpha-limit}
\lim\limits_{n\to\infty}\Delta_n=\alpha\in[0,\infty].
\end{equation}
Define $a_{n}=a(\Delta_n^{-1})$ and $b_{n}=b(\Delta_n^{-1})$ if $\alpha=0$, and $a_{n}=0$, $b_{n}=1$ if $\alpha\in\left( 0,\infty\right] $. Then
\begin{equation}\label{e17}
\lim\limits_{n\rightarrow\infty}P\Big(2\Delta_n^{-1/2}\{\max\limits_{1\leq j\leq n}\ln\left| \bz_{j}\right|-\frac12\sum\limits_{r=1}^{m_{n}}\psi(n_{r})\}\leq a_{n}+b_{n}y\Big)=\Phi_{\alpha}(y)
\end{equation}
for $y\in \mathbb{R}$.
\end{theorem}

Under condition \eqref{alpha-limit} with $\alpha\in [0,\infty)$, we have the limiting distribution for $\max\limits_{1\le j\le n}|\bz_j|$.

\begin{theorem}\label{t2} Assume condition \eqref{alpha-limit} hold with $\alpha\in [0,\infty)$.\\
(a). If $\alpha=0$, then $\alpha_{n}\big((\prod\limits_{r=1}^{m_{n}}n_{r})^{-1/2}\max\limits_{1\leq j\leq n}\left| \bz_{j}\right|-1\big)-\beta_{n}$ converges weakly to the Gumbel distribution $\Lambda(x)=\exp(-e^{-x})$, where
$\alpha_{n}=2\Delta_n^{-1/2}(-\ln\Delta_n)^{1/2}$ and $  \beta_{n}=-\ln\Delta_n-\ln(-\ln\Delta_n)-\ln\sqrt{2\pi}. $\\
(b). If $\alpha\in(0,\infty)$, then $(\prod\limits_{r=1}^{m_{n}}n_{r})^{-1/2}\max\limits_{1\leq j\leq n}\left| \bz_{j}\right|$ converges weakly to the $cdf$ $\Phi_{\alpha}(\alpha^{1/2}/2+2\alpha^{-1/2}\ln x)$, $x>0$.
\end{theorem}

\noindent{\bf Remark 1.} We can show under condition \eqref{alpha-limit} with $\alpha=
\infty$ that $\big(\max_{1\le j\le n}|\bz_j|-A_n\big)/B_n$ does not converge in distribution to any non-degenerate distribution for any  normalization constants $A_n\in \mathbb{R}$ and $B_n>0$.

\noindent{\bf Remark 2.} Under assumption $n=n_1=\cdots=n_{m_n+1}$,  the product ensemble $A_n^{(m_n)}$ is the product of $m_n$ independent Ginibre ensembles. In this case,  $\Delta_n=m_n/n$, and thus condition \eqref{alpha-limit} is equivalent to $\lim_{n\to\infty}m_n/n=\alpha\in [0,\infty]$. Then our Theorems~\ref{t1} and \ref{t2} reduce to, respectively, Proposition 2.1 and Theorem 3 in Jiang and Qi~\cite{JiangQi2017}.

\vspace{10pt}

Since $n_r\ge n$ for all $1\le r\le m_n$, we have $\Delta_n\le \sum^{m_n}_{r=1}1/n=m_n/n$. Hence  $\lim_{n\to\infty}m_n/n=0$ implies $\lim_{n\to\infty}\Delta_n=0$.  From Theorem~\ref{t2}, the limiting spectral radii is always Gumbel if $\lim_{n\to\infty}m_n/n=0$.  We have the following corollary.

\begin{coro}\label{cor1} Assume $\lim_{n\to\infty}m_n/n=0$. Then
$\alpha_{n}\big((\prod\limits_{r=1}^{m_{n}}n_{r})^{-1/2}\max\limits_{1\leq j\leq n}\left| \bz_{j}\right|-1\big)-\beta_{n}$ converges weakly to the Gumbel distribution $\Lambda(x)=\exp(-e^{-x})$, where
$\alpha_{n}=2\Delta_n^{-1/2}(-\ln\Delta_n)^{1/2}$ and $\beta_{n}=-\ln\Delta_n-\ln(-\ln\Delta_n)-\ln\sqrt{2\pi}$.
\end{coro}

To conclude this section, we provide some comments on the strategy for the proofs which are given in Section~\ref{proofs}.

\vspace{5pt}

\noindent{\bf Strategy for the proofs.} Much of our effort will be put in the proof of Theorem~\ref{t1}. We will first use a distributional representation for the spectral radii (see Lemmas~\ref{p1} below) and demonstrate that
the largest absolute eigenvalue has the same distribution as the maximum of $n$ products of independent Gamma random variables,  which implies that the logarithmic spectral radius has the same distribution as the maximum of sums of logarithmic Gamma random variables.
Then we decompose each sum of $m$ logarithmic Gamma random variables as a weighted sum of independent random variables plus a reminder term. Finally, we estimate the remainder (Lemmas~\ref{l5} and \ref{l7}) and apply
moderate deviation theorems to the weighted sums so as to estimate tail probabilities (see Lemmas~\ref{l9} and \ref{l10} below).
Somewhat similar steps here can be found in the proof of Proposition 2.1 in Jiang and Qi~\cite{JiangQi2017}, but our proofs are much more complicated as  we have to handle more parameters $n_1,\cdots, n_m$ other than only one parameter $m$ in Jiang and Qi~\cite{JiangQi2017}. For this reason we have to handle sum of weighted random variables in this paper (see, e.g.  Lemma~\ref{l10}) and employ new techniques to get finer estimates for remainders and tail probabilities (Lemmas~\ref{l7} and \ref{l8}).

\section{Proofs}\label{proofs}

In this section, we prove the main results given in Section~\ref{main}.   We first give some preliminary lemmas in Section~\ref{proof1}, and then provide the proofs for Theorems~\ref{t1} and \ref{t2} in Section~\ref{proof2}.

\subsection{Some Preliminary Lemmas}\label{proof1}

Define for $k>0$
\begin{equation}\label{delatjk}
\Delta_{j,k}=\sum\limits_{r=1}^{m_{n}}\frac{1}{(j+l_{r})^{k}}, ~~~ j=1,2,\cdots,n
\end{equation}
Note that
\[
\Delta_{n,k}=\sum^{m_n}_{r=1}\frac{1}{n_r^k}~~\mbox{ and }~~\Delta_n=\Delta_{n,1}.
\]

\begin{lemma} \label{p1}
Let $\left\{s_{j,r}, 1\leq r\leq m_{n},j\geq1\right\}$ be independent random variables and $s_{j,r}$ have the Gamma density $y^{j+l_{r}-1}e^{-y}I(y\geq0)/(j+l_{r}-1)!$ for each $j$ and $r$. Then $\max\limits_{1\leq j\leq n}{\left|\bz_{j}\right|^{2}}$ and $\max\limits_{1\leq j\leq n}{\prod\limits_{r=1}^{m_{n}}{s_{j,r}}}$ have the same distribution.
\end{lemma}

\noindent{\it Proof.} The lemma follows from Lemma 2.2 in Zeng~\cite{Zeng2017}. \hfill$\blacksquare$

\begin{lemma}\label{prod2sum} (Lemma 3.1 in Gui and Qi~\cite{GuiQi2018}) Suppose $\{l_n,\, n\geq 1\}$ is sequence of positive integers. Let $z_{nj}\in [0,1)$ be real numbers for $1\leq j \leq l_n$ such that $\max_{1\le j\le l_n}z_{nj}\to 0$ as $n\to\infty$. Then $\displaystyle\lim_{n\to\infty}\prod^{l_n}_{j=1}(1-z_{nj})\in (0,1)$ exists if and only if the limit $\displaystyle\lim_{n\to\infty}\sum^{l_n}_{j=1}z_{nj}=:z\in (0,\infty)$ exists and the relationship of the two limits is given by
\begin{equation}\label{prod-sum}
\lim_{n\to\infty}\prod^{l_n}_{i=1}(1-z_{ni})=e^{-z}.
\end{equation}
\end{lemma}

\begin{lemma}\label{l3} (Lemma 2.1 in Jiang and Qi~\cite{JiangQi2017}) Let $a_{ni}\in [0,1)$ be constants for $i\geq1$, $n\geq1$ and $\sup_{n\geq1,i\geq1}a_{ni}<1$. For each $i\geq1$, $a_{i}=\lim\limits_{n\to\infty}a_{ni}$. Assume $c_{n}=\sum\limits_{i=1}^{\infty}a_{ni}<\infty$ for each $n\geq1$ and $c=\sum\limits_{i=1}^{\infty}a_{i}<\infty$, and $\lim\limits_{n\to\infty}c_{n}=c$. Then,
\begin{equation*}
\lim\limits_{n\to\infty}\prod\limits_{i=1}^{\infty}(1-a_{ni})=\prod\limits_{i=1}^{\infty}(1-a_{i}).
\end{equation*}
\end{lemma}

\begin{lemma} \label{l4} (Lemma 2.2 in Jiang and Qi~\cite{JiangQi2017}) Let $\left\lbrace j_{n}, n\geq 1\right\rbrace $ and $\left\lbrace x_{n}, n\geq1\right\rbrace $ be positive numbers with $\lim\limits_{n\rightarrow\infty}x_{n}=\infty$ and $\lim\limits_{n\rightarrow\infty}j_{n}x_{n}^{-1/2}(\ln x_{n})^{1/2}=\infty$. For fixed $y\in\mathbb{R}$, if $\left\lbrace c_{n,j}, 1\leq j\leq j_{n}, n\geq1\right\rbrace $ are real numbers such that $\lim\limits_{n\rightarrow\infty}\max_{1\leq j\leq j_{n}}| c_{n,j}x_{n}^{1/2}-1| =0$, then
\begin{equation}\label{e4}
\lim\limits_{n\rightarrow\infty}\sum\limits_{j=1}^{j_{n}}(1-\Phi((j-1)c_{n,j}+a(x_{n})+b(x_{n})y))=e^{-y},\end{equation} where
$a(\cdot)$ and $b(\cdot)$ are defined in \eqref{ab}.
\end{lemma}

\begin{lemma} \label{l5} Set $G_{j}=\prod\limits_{r=1}^{m_{n}}{s_{j,r}}$, $1\leq j\leq n$, define the function $\eta(x)=x-1-\ln x$ for $x>0$, and write
\begin{equation}\label{e5}
M_{n}(i)=\max\limits_{n-i+1\leq j\leq n}{\Big|\sum\limits_{r=1}^{m_{n}}\big(\eta(\frac{s_{j,r}}{j+l_{r}})-E(\eta(\frac{s_{j,r}}{j+l_{r}}))\big)\Big|}.
\end{equation}
Recall $\psi(x)=\frac{\Gamma'{(x)}}{\Gamma(x)}$ as in \eqref{psi}.  Then for $1\le i\le n$
\[
\Big|\max\limits_{n-i+1\leq j\leq n}{\ln G_{j}}-\max\limits_{n-i+1\leq j\leq n}\Big(\sum\limits_{r=1}^{m_{n}}{\frac{s_{j,r}-(j+l_{r})}{j+l_{r}}}+\sum\limits_{r=1}^{m_{n}}{\psi(j+l_{r})}\Big)\Big|\leq M_{n}(i).
\]
\end{lemma}

\noindent{\it Proof}.
The moment-generating function of $\ln s_{j,r}$ is
\begin{equation}\label{e6}
m_{j,r}=E(e^{t\ln s_{j,r}})=\frac{\Gamma(j+l_{r}+t)}{\Gamma(j+l_{r})}
\end{equation} for $t>-j-l_{r}$.
Then, we have
\begin{equation}\label{e7}
E(\ln s_{j,r})=\frac{d}{dt}m_{j,r}(t)|_{t=0}=\frac{\Gamma'(j+l_{r})}{\Gamma(j+l_{r})}=\psi(j+l_{r}).
\end{equation}
Using the relationship  $\ln x=x-1-\eta(x)$, we can rewrite $\ln G_{j}$ as
\begin{align*}
\ln G_j&=\ln\prod_{r=1}^{m_{n}}s_{j,r}\\
&=\sum^{m_n}\limits_{r=1}\ln\frac{s_{j,r}}{j+l_{r}}
+\sum\limits_{r=1}^{m_{n}}\ln(j+l_{r})\\&=\sum\limits_{r=1}^{m_{n}}\frac{s_{j,r}-(j+l_{r})}{j+l_{r}}
-\sum\limits_{r=1}^{m_{n}}\eta(\frac{s_{j,r}}{j+l_{r}})+\sum\limits_{r=1}^{m_{n}}\ln(j+l_{r})\\
&=\sum\limits_{r=1}^{m_{n}}\frac{s_{j,r}-(j+l_{r})}{j+l_{r}}+\sum\limits_{r=1}^{m_{n}}\psi(j+l_{r})
-\sum\limits_{r=1}^{m_{n}}\Big(\eta(\frac{s_{j+l_{r}}}{j+l_{r}})-\ln(j+l_{r})+\psi(j+l_{r})\Big).
\end{align*}
Since $E(\ln s_{j,r})=\psi(j+l_{r})$ from \eqref{e7}, we obtain that
\begin{equation}\label{etamean}
E(\eta(\frac{s_{j,r}}{j+l_{r}}))=\ln(j+l_{r})-\psi(j+l_{r}),
\end{equation}
and thus we have,
\begin{equation}\label{e8}
 \ln G_{j}= \sum\limits_{r=1}^{m_{n}}\frac{s_{j,r}-(j+l_{r})}{j+l_{r}}+
\sum\limits_{r=1}^{m_{n}}\psi(j+l_{r})-\sum\limits_{r=1}^{m_{n}}\Big(\eta(\frac{s_{j,r}}{j+l_{r}})-E(\eta(\frac{s_{j,r}}{j+l_{r}}))\Big).
\end{equation}
Note that for any two sequences of real numbers
$\left\lbrace x_{n}\right\rbrace $ and $\left\lbrace y_{n}\right\rbrace $,
\[
\big| \max\limits_{1\leq j\leq n}x_{j}-\max\limits_{1\leq j\leq n}y_{j}\big|\leq\max\limits_{1\leq j\leq n}\big| x_{j}-y_{j}\big|.
\]
Then it follows from (\ref{e8}) that
 \[
\big|\max\limits_{n-i+1\leq j\leq n}{\ln G_{j}}-\max\limits_{n-i+1\leq j\leq n}\big(\sum\limits_{r=1}^{m_{n}}{\frac{s_{j,r}-(j+l_{r})}{j+l_{r}}}+\sum\limits_{r=1}^{m_{n}}{\psi(j+l_{r})}\big)\big|\leq M_{n}(i).\]
This complete the proof of the lemma. \hfill$\blacksquare$

\begin{lemma} \label{l6} Recall $\Delta_{n,j}$ is defined in \eqref{delatjk}. Assume $\left\{j_{n};n\geq1\right\}$ is a sequence of numbers satisfying $1\leq j_{n}\leq n/2$ for all $n\ge 2$, then for $n-j_{n}+1\leq j\leq n$, we have \\
(1) $\Delta_{n,k}\leq\Delta_{j,k}<2^{k}\Delta_{n,k}$ for any $k>0$;\\
(2) $\Delta_{j,2}/\Delta_{j,1}^{1+a}\leq j^{a-1}$ for any $a\ge 0$.
\end{lemma}

\noindent{\it Proof}. Assume $n-j_{n}+1\leq j\leq n$.  Since $\frac{n_{r}}{2}<n_{r}-j_{n}+1\leq j+l_{r}\leq n_{r}$, we have for $k>0$,
\[
\frac{1}{n_{r}^{k}}\leq\frac{1}{(j+l_{r})^{k}}<\frac{2^{k}}{n_{r}^{k}}, ~~~~~1\le r\le m_n.
\]
By summing up over $r\in \{1, \cdots, m_n\}$, we obtain that $\Delta_{n,k}\leq\Delta_{j,k}<2^{k}\Delta_{n,k}$, i.e.  (1) holds.

Note that $l_r\ge 0$ and $l_1=0$. We have that $j/(j+l_r)\le 1$ for any $1\le j\le n$ and $1\le r\le m_n$, and $\Delta_{n,j}\ge 1/j$. Therefore,
for any $a\ge 0$,
\[
\frac{\Delta_{j,2}}{\Delta_{j,1}^{1+a}}=\frac{\sum\limits_{r=1}^{m_{n}}\frac{1}{(j+l_{r})^{2}}}{(\sum\limits_{r=1}^{m}\frac{1}{j+l_{r}})^{1+a}}
=j^{a-1}\cdot\frac{\sum\limits_{r=1}^{m_{n}}(\frac{j}{j+l_{r}})^{2}}{(\sum\limits_{r=1}^{m_{n}}\frac{j}{j+l_{r}})^{1+a}}
\le j^{a-1}\cdot\frac{\sum\limits_{r=1}^{m_{n}}\frac{j}{j+l_{r}}}{(\sum\limits_{r=1}^{m_{n}}\frac{j}{j+l_{r}})^{1+a}}\le \frac{j^{a-1}}
{(\sum\limits_{r=1}^{m_{n}}\frac{j}{j+l_{r}})^{a}} \le j^{a-1}.
\]
In the last estimation we have used the fact that $\sum\limits_{r=1}^{m_{n}}\frac{j}{j+l_{r}}\ge \frac{j}{j+l_1}=1$.\hfill$\blacksquare$

\begin{lemma} \label{l7} Assume $\left\{j_{n}, ~n\geq1\right\}$ is a sequence of numbers satisfying $1\leq j_{n}\leq n/2$ for all $n\ge 2$. Then, $M_{n}(j_{n})=O_{p}(j_{n}(\frac{\Delta_n}{n})^{1/2})$ and  $M_{n}(j_{n})=O_{p}(\Delta_n\ln n)$ as $n\to\infty$.
\end{lemma}

\noindent{\it Proof}. We have
\begin{align*}
E(M_{n}(j_{n}))&\leq\sum\limits_{j=n-j_{n}+1}^{n}E\Big|\sum\limits_{r=1}^{m_{n}}
\big(\eta(\frac{s_{j,r}}{j+l_{r}})-E(\eta(\frac{s_{j,r}}{j+l_{r}}))\big)\Big|\\
&\leq\sum\limits_{j=n-j_{n}+1}^{n}\Big\{E\Big(\sum\limits_{r=1}^{m_{n}}
\big(\eta(\frac{s_{j,r}}{j+l_{r}})-E(\eta(\frac{s_{j,r}}{j+l_{r}}))\big)\Big)^2\Big\}^{1/2}\\
&=\sum\limits_{j=n-j_{n}+1}^{n}\Big\{\sum\limits_{r=1}^{m_{n}}E\big(\eta(\frac{s_{j,r}}{j+l_{r}})-E(\eta(\frac{s_{j,r}}{j+l_{r}}))\big)^{2}\Big\}^{1/2}\\
&\leq\sum\limits_{j=n-j_{n}+1}^{n}\Big\{\sum\limits_{r=1}^{m_{n}}E\big(\eta(\frac{s_{j,r}}{j+l_{r}})\big)^{2}\Big\}^{1/2}\\
&\leq\sum\limits_{j=n-j_{n}+1}^{n}\Big\{\sum\limits_{r=1}^{m_{n}}\frac{\big(\frac{s_{j,r}}{j+l_{r}}-1\big)^{4}}{(2\min(\frac{s_{j,r}}{j+l_{r}},1))^{2}}
\Big\}^{1/2}\\
&=\frac12\sum\limits_{j=n-j_{n}+1}^{n}\Big\{\sum\limits_{r=1}^{m_{n}}E((\frac{s_{j,r}-(j+l_{r})}{j+l_{r}})^{4})(\min(\frac{s_{j,r}}{j+l_{r}},1))^{-2}
\Big\}^{1/2}.
 \end{align*}
In the last inequality we have used estimation that
\[
0\le \eta(x)=x-1-\ln x=\int^x_1\frac{t-1}{t}dt\le \frac{(x-1)^2}{2\min(x,1)},  ~~~x>0.
\]
Since $s_{j,r}$ has density $y^{j+l_{r}-1}e^{-y}I(y>0)/(j+l_{r}-1)!$, we have $E(s_{j,r}^{-4})=\frac{\Gamma(j+l_{r}-4)}{\Gamma(j+l_{r})}$. By
the Marcinkiewicz-Zygmund inequality(see, for example, Corollary 2 in Section 10.3 from Chow and Teicher~\cite{chow2012}), we obtain $E(s_{j,r}-(j+l_{r}))^{8}\leq C(j+l_{r})^{4}$, where $C$ is a constant not depending on $j$. From now on we will use $C$ to denote a generic constant which may be different at different places.  Then we have
\begin{align*}
&E((\frac{s_{j,r}-(j+l_{r})}{j+l_{r}})^{4}(\min(\frac{s_{j,r}}{j+l_{r}},1))^{-2})\\
\allowdisplaybreaks[4]
&\leq(E(\frac{s_{j,r}-(j+l_{r})}{j+l_{r}})^{8}\cdot E(\min(\frac{s_{j,r}}{j+l_{r}},1))^{-4})^{1/2}\\
&\leq(E(\frac{s_{j,r}-(j+l_{r})}{j+l_{r}})^{8}\cdot E\big(1+(\frac{j+l_{r}}{s_{j,r}})^{4})\big)^{1/2}\\
&\leq\Big(1+\frac{(j+l_{r})^{3}}{(j+l_{r}-1)(j+l_{r}-2)(j+l_{r}-3)}\Big)^{1/2}\Big(E\big(\frac{s_{j,r}-(j+l_{r})}{j+l_{r}}\big)^{8}\Big)^{1/2}\\
&\leq C(j+l_{r})^{-2},
\end{align*}
and thus from Lemma~\ref{l6} we obtain
\begin{align*}
E(M_{n}(j_{n}))&\leq\frac{\sqrt{C}}{2}\sum\limits_{j=n-j_{n}+1}^{n}\big(\sum\limits_{r=1}^{m_{n}}(j+l_{r})^{-2}\big)^{1/2}\\
&=\frac{\sqrt{C}}{2}\sum\limits_{j=n-j_{n}+1}^{n}\Delta_{j,2}^{1/2}\\
&\le\sqrt{C}\sum\limits_{j=n-j_{n}+1}^{n}\Delta_{n,2}^{1/2}\\
&=\sqrt{C}j_n\Delta_{n,2}^{1/2}\\
&\le O(\frac{j_n}{n^{1/2}}\Delta_{n,1}^{1/2}).
\end{align*}
Therefore $M_{n}(j_{n})=O_{p}\big(\frac{j_n}{n^{1/2}}\Delta_{n}^{1/2}\big)$.\\

Recall $\psi(x)=\frac{\Gamma'(x)}{\Gamma(x)}$ for $x>0$. By Formulas 6.3.18 and 6.4.12 in Abramowitz and Stegun~\cite{Abramowitz1972} we have
\begin{equation}\label{e9}\psi(x)=\ln x-\frac{1}{2x}+O(\frac{1}{x^{2}})\ and\ \psi'(x)=\frac{1}{x}+\frac{1}{2x^{2}}+O(\frac{1}{x^{3}})\end{equation}
as $x\rightarrow +\infty$. From \eqref{etamean}, $E\eta(\frac{s_{j,r}}{j+l_{r}})=\ln(j+l_{r})-\psi(j+l_{r})=O(\frac{1}{j+l_{r}})$ as $j\rightarrow\infty$, we have
\begin{equation}\label{e10}
M_{n}(j_{n})\leq\max\limits_{n-j_{n}+1\leq j\leq n}\sum\limits_{r=1}^{m_{n}}\eta(\frac{s_{j,r}}{j+l_{r}})+O(\sum\limits_{r=1}^{m_{n}}\frac{1}{n_{r}}).
\end{equation}

For $n-j_{n}+1\leq j\leq n$, we consider the moment generating function of $\eta(\frac{s_{j,r}}{j+l_{r}})$. Since $s_{j,r}$ has a Gamma($j+l_{r}$) distribution, we have
\begin{align*} Ee^{t\eta(\frac{s_{j,r}}{j+l_{r}})}&=E\Big(\exp\big(t(\frac{s_{j,r}}{j+l_{r}}-1-\ln\frac{s_{j,r}}{j+l_{r}})\big)\Big)\\
&=e^{-t}E\big((\frac{s_{j,r}}{j+l_{r}})^{-t}\exp(t\cdot\frac{s_{j,r}}{j+l_{r}})\big)\\
&=\frac{e^{-t}(j+l_{r})^{t}}{\Gamma(j+l_{r})}\int_{0}^{\infty}x^{j+l_{r}-t-1}e^{-x(1-\frac{t}{j+l_{r}})}dx\\
&=\frac{e^{-t}(j+l_{r})^{t}}{\Gamma(j+l_{r})}\int_{0}^{\infty}(\frac{j+l_{r}}{j+l_{r}-t})^{j+l_{r}-t}y^{j+l_{r}-t-1}e^{-y}dy\\
&=e^{-t}(j+l_{r})^{t}\frac{\Gamma(j+l_{r}-t)}{\Gamma(j+l_{r})}(\frac{j+l_{r}}{j+l_{r}-t})^{j+l_{r}-t}.
\end{align*}
Uniformly over $0<t<n/4$, we have from \eqref{e9}
\begin{align*}
\ln\frac{\Gamma(j+l_{r}-t)}{\Gamma(j+l_{r})}
&=\int_{j+l_{r}}^{j+l_{r}-t}\psi(x)dx=\int_{j+l_{r}}^{j+l_{r}-t}(\ln x-\frac{1}{2x}+O(\frac{1}{x^{2}}))dx\\
&=(x\ln x-x)|_{j+l_{r}}^{j+l_{r}-t}-\frac12\ln\frac{j+l_{r}-t}{j+l_{r}}+O(\frac{t}{(j+l_{r}-t)^{2}})\\
&=(j+l_{r}-t)\ln(j+l_{r}-t)-(j+l_{r})\ln(j+l_{r})+t\\&-\frac12\ln\frac{j+l_{r}-t}{j+l_{r}}+O(\frac{t}{(j+l_{r})^{2}}).
\end{align*}
Therefore, we obtain
\begin{equation}\label{e11}
\frac{\Gamma(j+l_{r}-t)}{\Gamma(j+l_{r})}
=e^{t}\frac{(j+l_{r}-t)^{j+l_{r}-t}}{(j+l_r)^{j+l_{r}}}(1-\frac{t}{j+l_{r}})^{-1/2}\exp\big(O(\frac{t}{(j+l_{r})^{2}})\big)
\end{equation}
and
\begin{align*}
E\exp\big(t\eta(\frac{s_{j,r}}{j+l_{r}})\big)&=(1-\frac{t}{j+l_{r}})^{-1/2}\exp\Big(O(\frac{t}{(j+l_{r})^{2}})\Big)\\
&=\exp\Big(\frac12\cdot\frac{t}{j+l_{r}}+\frac{1}{4}\cdot\frac{t^{2}}{(j+l_{r})^{2}}+O(\frac{t}{(j+l_{r})^{2}})\Big).
\end{align*}
Then we have
\begin{align*}
E\exp\big(t\sum\limits_{r=1}^{m_{n}}\eta(\frac{s_{j,r}}{j+l_{r}})\big)&=\exp\big(\frac{t}2\Delta_{j,1}+O(\Delta_{j,2}t+\Delta_{j,2}t^{2})\big)\\
&\leq\exp\big(t\Delta_n+O(\Delta_{n,2}t^{2}+\Delta_{n,2}t)\big)
\end{align*}
uniformly over $0<t<n/4$ and $n-j_n+1\le j\le n$ as $n\to\infty$.
Now plug in $t=1/(4\Delta_n)$. Since  $\Delta_n\geq\frac{1}{n}$, we have  $0<t\leq\frac{n}{4}$, and thus we get
\begin{align*}
&P(\sum\limits_{r=1}^{m_{n}}\eta(\frac{s_{j,r}}{j+l_{r}})>8\Delta_n\ln n)\\
&\leq\frac{E\big(\exp(t\sum\limits_{r=1}^{m_{n}}\eta(\frac{s_{j,r}}{j+l_{r}}))\big)}{\exp(8t\Delta_n\ln n)}\\
&\leq\frac{\exp\big(4+O(\Delta_{n,2}/\Delta_{n,1}^2+\Delta_{n,2}/\Delta_{n,1})\big)}{\exp(2\ln n)}\\
&=O(n^{-2})
\end{align*}
from Lemma~\ref{l6}. Therefore,
$$P\Big(\max\limits_{n-j_{n}+1\leq j\leq n}\sum\limits_{r=1}^{m_{n}}\eta(\frac{s_{j,r}}{j+l_{r}})>8\Delta_n\ln n\Big)\leq O(n^{-1})\to 0, $$
which means
$$M_{n}(j_{n})\leq\max\limits_{n-j_{n}+1\leq j\leq n}\sum\limits_{r=1}^{m_{n}}\eta(\frac{s_{j,r}}{j+l_{r}})+O(\Delta_n)=O_{p}(\Delta_n\ln n).$$
This completes the proof. \hfill$\blacksquare$

\begin{lemma} \label{l8}
Let $\left\lbrace j_{n}, ~ n\geq 1\right\rbrace $ be positive integers satisfying
\begin{equation}\label{assup1}
\lim\limits_{n\rightarrow\infty}{\frac{j_{n}}{n}}=0, ~~\lim\limits_{n\rightarrow\infty}j_{n}(\frac{\Delta_n}{\ln n})^{1/2}=\infty.
\end{equation}
Then, for any $x\in \mathbb{R}$
\begin{equation}\label{e12}
\lim\limits_{n\rightarrow\infty}\sum\limits_{j=1}^{n-j_{n}}P(\ln G_{j}>\sum\limits_{r=1}^{m_{n}}\psi(n+l_{r})+\Delta_n^{1/2}x)=0.
\end{equation}
\end{lemma}

\noindent{\it Proof}. Fix $x\in\mathbb{R}$. For each $1\leq j\leq n-j_{n}$ and any $t>0$, we have from \eqref{e6} that
\begin{eqnarray*}
&&P(\ln G_{j}>\sum\limits_{r=1}^{m_{n}}\psi(n+l_{r})+\Delta_n^{1/2}x)\\
&\leq&\frac{E(e^{t\ln G_{j}})}{\exp(t(\sum\limits_{r=1}^{m_{n}}\psi(n+l_{r})+\Delta_n^{1/2}x))}\\
&=&\exp\Big(\sum\limits_{r=1}^{m_{n}}(\ln\Gamma(j+l_{r}+t)-\ln\Gamma(j+l_{r}))-t(\sum\limits_{r=1}^{m_{n}}\psi(n+l_{r})
+\Delta_n^{1/2}x)\Big)\\
&=&\exp\Big(\sum\limits_{r=1}^{m_{n}}\int_{0}^{t}(\psi(j+l_{r}+s)-\psi(j+l_{r}))ds-t(\sum\limits_{r=1}^{m_{n}}(\psi(n+l_{r})
-\psi(j+l_{r}))+\Delta_n^{1/2}x)\Big).
\end{eqnarray*}

Since there exists an integer $j_{0}$ such that for all $j_{0}\leq j\leq n-j_{n}$ and for all $1\leq r\leq m_{n}$, $$\ln\frac{j+l_{r}+s}{j+l_{r}}\leq\psi(j+l_{r}+s)-\psi(j+l_{r})=\int_{0}^{s}\psi'(j+l_{r}+v)dv\leq\frac{1.1s}{j+l_{r}}.$$
By the first inequality above, for all  $j_{0}\leq j\leq n-j_{n}$, $1\le r\le m_n$ and all large $n$, $$\psi(n+l_{r})-\psi(j+l_{r})\geq\ln\frac{n+l_{r}}{j+l_{r}}\geq\ln\frac{n_{r}}{n_{r}-j_{n}}=-\ln(1-\frac{j_{n}}{n_{r}})\geq\frac{0.999j_{n}}{n_{r}},$$
which implies
\[
\sum\limits_{r=1}^{m_{n}}(\psi(n+l_{r})-\psi(j+l_{r}))\geq \sum^{m_n}_{r=1}\ln\frac{n_r}{j+l_{r}}
\]
and
\[
\sum\limits_{r=1}^{m_{n}}(\psi(n+l_{r})-\psi(j+l_{r}))\geq 0.999j_n\Delta_n
\]
uniformly for $j_{0}\leq j\leq n-j_{n}$ for all large $n$.
By assumption \eqref{assup1}, we have $\Delta_n^{1/2}=o(j_{n}\Delta_n)$, and
 $$
\sum\limits_{r=1}^{m_{n}}(\psi(n+l_{r})-\psi(j+l_{r}))+\Delta_n^{1/2}x\geq0.99\sum\limits_{r=1}^{m_{n}}\ln\frac{n_{r}}{j+l_{r}}$$ uniformly over $j_{0}\leq j\leq n-j_{n}$ for all large $n$. Therefore, for all $j_{0}\leq j\leq n-j_{n}$,
\begin{eqnarray*}
&&P(\ln G_{j}>\sum\limits_{r=1}^{m_{n}}\psi(n+l_{r})+\Delta_n^{1/2}x)\\
&\leq&\exp\Big\{ 1.1\sum\limits_{r=1}^{m_{n}}\int_{0}^{t}\frac{s}{j+l_{r}}ds-0.99t\sum\limits_{r=1}^{m_{n}}\ln \frac{n_{r}}{j+l_{r}}\Big\}\\
&=&\exp\Big\{\sum\limits_{r=1}^{m_{n}}\frac{0.55t^{2}}{j+l_{r}}-0.99t\sum\limits_{r=1}^{m_{n}}\ln \frac{n_{r}}{j+l_{r}}\Big\} \\
&=&\exp\Big\{0.55t^{2}\Delta_{j,1}-0.99t\sum\limits_{r=1}^{m_{n}}\ln \frac{n_{r}}{j+l_{r}}\Big\}
\end{eqnarray*}
for all $t>0$ and large $n$. By selecting $t=0.9\sum\limits_{r=1}^{m_{n}}\ln\frac{n_{r}}{j+l_{r}}/\Delta_{j,1}$, we have
\begin{equation}\label{e13}
P(\ln G_{j}>\sum\limits_{r=1}^{m_{n}}\psi(n+l_{r})+\Delta_n^{1/2}x)
\leq\exp\Big\{-\frac{0.4455}{\Delta_{j,1}}\Big(\sum\limits_{r=1}^{m_{n}}\ln\frac{n_{r}}{j+l_{r}}\Big)^{2}\Big\}
\end{equation}
uniformly over $j_{0}\leq j\leq n-j_{n}$ for all large $n$.

Now we turn to estimate the probability on the right-hand side  of \eqref{e13}. For each $r\in \{1, \cdots, m_n\}$, define the function
$f_r(x)=x(\ln n_r-\ln x)$, $0< x\le n_r$. Note that $f_r'(x)=\ln n_r-\ln x-1$ is decreasing and $f_r''(x)=-1/x<0$
for $x\in (0,  n_r]$.  This implies that $f_r(x)$ is concave in $x\in (0,  n_r]$, and for any constants $0<a<b<n_r$, the minimum value of $f_r(x)$ over $[a, b]$
is achieved at the two endpoints of interval $[a,b]$, i.e.,
\begin{equation}\label{minimal}
\min_{a\le x\le b}f_r(x)=\min\big(f_r(a), f_r(b)\big).
\end{equation}

For any $1\le j\le n-j_n$ and $1\le r\le m_n$, set $a_{nj}=\min(j, n/8)$ and $b_{nj}=n_r-j_n$. Then  $1\le a_{nj}\le j+l_r\le b_{nj}<n_r$ holds uniformly
over $1\le j\le n-j_n$ and $1\le r\le m_n$ for for all large $n$.
Note that
\[
f_r(a_{nj})=a_{nj}\ln\frac{n_r}{a_{nj}}\ge a_{nj}\ln\frac{n}{a_{nj}}
\]
and
\[
f_r(b_{nj})\ge (n-j_n)\ln\frac{n_r}{n_r-j_n}= -(n-j_n)\ln(1-\frac{j_n}{n_r})\ge-(n-j_n)\ln(1-\frac{j_n}{n}) \ge \frac12j_n
\]
for all large $n$.  By applying
\eqref{minimal} we obtain from \eqref{minimal} that
\[
(j+l_r)\ln\frac{n_r}{j+l_r}\ge \min(a_{nj}\ln\frac{n}{a_{nj}}, \frac{j_n}{2})=:\delta_{nj},
\]
or equivalently
\[
\ln\frac{n_r}{j+l_r}\ge\frac{\delta_{nj}}{j+l_r}
\]
over $1\le j\le n-j_n$ and $1\le r\le r\le m_n$ for all large $n$.  Therefore, we conclude that
\begin{equation}\label{allj}
\sum\limits_{r=1}^{m_{n}}\ln\frac{n_{r}}{j+l_{r}}\ge \delta_{nj}\sum\limits_{r=1}^{m_{n}}\frac{1}{j+l_{r}}=\delta_{nj}\Delta_{j,1}
\end{equation}
uniformly over $1\le j\le n-j_n$ for all large $n$.  Thus, for all large $n$,
\begin{eqnarray}\label{minmin}
\min\limits_{1\leq j\leq n-j_{n}}\Delta_{j,1}^{-1}\Big(\sum\limits_{r=1}^{m_{n}}\ln \frac{n_{r}}{j+l_{r}}\Big)^{2}
&\ge&\min_{1\le j\le n-j_n}\delta_{nj}^2\Delta_{j,1}\nonumber\\
&=&\min_{1\le j\le n-j_n}\min\big(a_{nj}^2(\ln\frac{n}{a_{nj}})^2\Delta_{j,1}, \frac14j_n^2\Delta_{j,1}\big)\nonumber\\
&\ge&\min_{1\le j\le n-j_n}\min\big(\frac18a_{nj}(\ln\frac{n}{a_{nj}})^2, \frac14j_n^2\Delta_n\big)\nonumber\\
&=&\min\big(\min_{1\le j\le n-j_n}\frac18a_{nj}(\ln\frac{n}{a_{nj}})^2, \frac14j_n^2\Delta_n\big).
\end{eqnarray}
To obtain the second inequality above we have used the facts that $\Delta_{j,1}\ge 1/j$, $a_{nj}/j=\min(j, n/8)/j\ge 1/8$ and $\Delta_{j,1}\ge \Delta_{n,1}=\Delta_n$.

Our aim is to show that
\begin{equation}\label{oederln(n)}
\frac{1}{\ln n}\min\limits_{1\leq j\leq n-j_{n}}\Delta_{j,1}^{-1}\Big(\sum\limits_{r=1}^{m_{n}}\ln \frac{n_{r}}{j+l_{r}}\Big)^{2}\to
\infty ~~\mbox{ as }n\to\infty.
\end{equation}
In fact, condition \eqref{assup1} implies $j_n^2\Delta_n/\ln n\to \infty$ as $n\to\infty$. By \eqref{minmin} it remains to show that
\begin{equation}\label{orderln(n)2}
\frac1{\ln n}\min_{1\le j\le n-j_n}a_{nj}(\ln\frac{n}{a_{nj}})^2\to\infty~~\mbox{ as }n\to\infty.
\end{equation}
To show this, we consider the function $f(x)=x(\ln n-\ln x)^2$, $1\le x\le n/8$. $f(x)$ is increasing since $f'(x)=(\ln n-\ln x))(\ln n-\ln x-2)>0$ for $x\in [0, n/8]$. Therefore, we have $\min_{1\le x\le n/8}f(x)\ge f(1)=(\ln n)^2$, which implies that $a_{nj}(\ln\frac{n}{a_{nj}})^2\ge (\ln n)^2$, and the left-hand side of \eqref{orderln(n)2} is larger than $\ln n$. This proves \eqref{orderln(n)2}.

Now it follows from \eqref{oederln(n)} that
\[
\min\limits_{j_{0}\leq j\leq n-j_{n}}\Delta_{j,1}^{-1}\Big(\sum\limits_{r=1}^{m_{n}}\ln \frac{n_{r}}{j+l_{r}}\Big)^{2}\ge 10\ln n
\]
for all large $n$, which coupled with \eqref{e13} implies
\[
\max_{j_{0}\leq j\leq n-j_{n}}P\big(\ln G_{j}>\sum\limits_{r=1}^{m_{n}}\psi(n+l_{r})+\Delta_n^{1/2}x\big)\leq\exp(-4.4\ln n)=n^{-4.4},
\]
and hence,
\[
\sum\limits_{j=j_{0}}^{n-j_{n}}P\big(\ln G_{j}>\sum\limits_{r=1}^{m_{n}}\psi(n+l_{r})+\Delta_n^{1/2}x\big)=O(n^{-3.4})\rightarrow 0
~~\mbox{ as }n\to\infty.
\]

Finally, we will consider the tail probability of $\ln G_j$ when $1\le j<j_0$. From \eqref{e6} we have
\[
E(G_j)=\prod^{m_n}_{r=1}\frac{\Gamma(j+l_r+1)}{\Gamma(j+l_r)}=\prod^{m_n}_{r=1}(j+l_r).
\]
Using \eqref{e9} we get for all large $n$
\begin{eqnarray*}
\sum\limits_{r=1}^{m_{n}}\psi(n+l_{r})+\Delta_n^{1/2}x&=&\sum^{m_n}_{r=1}\ln(n+l_r)+O(\Delta_n+\Delta_n^{1/2})\\
&\ge& \sum^{m_n}_{r=1}\ln(n+l_r)+O(\Delta_n+1).
\end{eqnarray*}
For each fixed $j$, $1\le j<j_0$, since $G_j>0$, we have from Chebyshev's inequality and equation \eqref{allj} that
\begin{eqnarray*}
&&P\big(\ln G_{j}>\sum\limits_{r=1}^{m_{n}}\psi(n+l_{r})+\Delta_n^{1/2}x\big)\\
&=&P\big(G_{j}>\exp\{\sum\limits_{r=1}^{m_{n}}\psi(n+l_{r})+\Delta_n^{1/2}x\}\big)\\
&\le &\frac{E(G_j)}{\exp\{\sum\limits_{r=1}^{m_{n}}\psi(n+l_{r})+\Delta_n^{1/2}x\}}\\
&\le &\exp\{-\sum\limits_{r=1}^{m_{n}}\ln\frac{n+l_{r}}{j+l_r}+O(\Delta_n+1)\}\\
&\le &\exp\{-(1+o(1))\sum\limits_{r=1}^{m_{n}}\ln\frac{n+l_{r}}{j+l_r}+O(1)\}\\
&\le &\exp\{-(1+o(1))\ln\frac{n+l_1}{j+l_1}+O(1)\}\\
&\le &\exp\{-(1+o(1))\ln\frac{n}{j}+O(1)\}\\
&\to& 0
\end{eqnarray*}
as $n\to\infty$.  This proves \eqref{e12} and completes the proof of the lemma. \hfill$\blacksquare$

\begin{lemma} \label{l9} (Proposition 4.5 in Chen, Fang and Shao~\cite{chen2013}) Let $\xi_{i}$, $1\leq i\leq n$ be independent random variables with $E\xi_{i}=0$ and $Ee^{t_{n}\left|\xi_{i}\right|}<\infty$,  $1\leq i\leq n$ for some $t_{n}$. Assume that $\sum\limits_{i=1}^{n}{E\xi_{i}^{2}}=1$. Then
\begin{equation}\label{e14}
\frac{P(W\geq x)}{1-\Phi(x)}=1+O(1)(1+x^{3})\gamma e^{4x^{3}\gamma}
\end{equation} for $0\leq x\leq t_{n}$, where $W=\sum\limits_{i=1}^{n}{\xi_{i}}$ and $\gamma=\sum\limits_{i=1}^{n}E\big(|\xi_{i}|^{3}e^{x|\xi_{i}|}\big)$.
\end{lemma}

\begin{lemma} \label{l10} Let $\left\lbrace j_{n}, n\geq 1\right\rbrace$ be positive integers satisfying $1\leq j_{n}\leq n/2$ and $\lim\limits_{n\rightarrow\infty}{\frac{j_{n}}{n}=0}$. Let $W_{j}=\Delta_{j,1}^{-1/2}\sum\limits_{r=1}^{m_{n}}\big(s_{j,r}-(j+l_{r})\big)/(j+l_r)$ and $t_n=O(n^{1/7})$ be any sequence of positive numbers.
 Then $P(W_{j}\geq x)=\big(1-\Phi(x)\big)(1+o(1))$ uniformly over $0\leq x\leq t_n$ and $n-j_{n}+1\leq j\leq n$
as $n\rightarrow\infty$.
\end{lemma}

\noindent{\it Proof}. Let $\{X_{i,r}, ~i\ge 1, r\ge 1\}$ be an array of i.i.d. random variables with the standard exponential distribution. Then for each $j$,
$\{s_{j,r},~ 1\le r\le m_n\}$ have the same joint distribution as $\{\sum^j_{i=1}X_{i,r}, 1\le r\le m_n\}$. Without loss of generality we assume $s_{j,r}=\sum^j_{i=1}X_{i,r}$ for $1\le r\le m_n$, $n-j_n\le j\le n$.

Set $d_{j,r}=(j+l_{r})^{-1}$ and $D_{j,r}=d_{j,r}/\Delta_{j,1}^{1/2}$ for $1\leq r\leq m_{n}$. Then
\begin{align*} W_{j}&=\Delta_{j,1}^{-1/2}\sum\limits_{r=1}^{m_{n}}\sum\limits_{i=1}^{j+l_{r}}{\frac{1}{(j+l_{r})}(X_{i,r}-1)}\\
&=\Delta_{j,1}^{-1/2}\sum\limits_{r=1}^{m_{n}}\sum\limits_{i=1}^{j+l_{r}}d_{j,r}(X_{i,r}-1)\\&=\sum\limits_{r=1}^{m_{n}}\sum\limits_{i=1}^{j+l_{r}}\xi_{i,r},
\end{align*}
where $\xi_{i,r}=D_{j,r}(X_{i,r}-1)$. Since $E(X_{i,r})=Var(X_{i,r})=1$, we obtain
\[
E\xi_{i,r}=0~~ \mbox{ and  } ~\sum\limits_{r=1}^{m_{n}}\sum\limits_{i=1}^{j+l_{r}}E\xi_{i,r}^{2}=1.
\]
 Furthermore, we have
\begin{align*}
\sum\limits_{r=1}^{m_{n}}\sum\limits_{i=1}^{j+l_{r}}E\big(\left|\xi_{1,r}\right|^{3}e^{t\left|\xi_{i,r}\right|}\big)
&=\sum\limits_{r=1}^{m_{n}}E(\left|D_{j,r}(X_{1,r}-1)\right|^{3}e^{t\left|D_{j,r}(X_{1,r}-1)\right|})\cdot\frac{1}{d_{j,r}}\\
&=\Delta_{j,1}^{-3/2}\sum\limits_{r=1}^{m_{n}}E(d_{j,r}^{3}\left|X_{1,r}-1\right|^{3}e^{tD_{j,r}\left|X_{1,r}-1\right|})\cdot\frac{1}{d_{j,r}}\\
&\leq\Delta_{j,1}^{-3/2}\sum\limits_{r=1}^{m_{n}}d_{j,r}^{2}E\big((X_{1,r}^{3}+1)(e^{tD_{j,r}(X_{1,r}-1)}+e^{-tD_{j,r}(X_{1,r}-1)})\big).
\end{align*}
Using the moment-generating function $E(e^{tD_{j,r}\cdot X_{i,r}})=(1-D_{j,r}t)^{-1}$, we have
 \[
E(X_{1,r}^{3}e^{tD_{j,r}\cdot X_{1,r}})=\frac6{(1-D_{j,r}t)^{4}},
\]
thus
\begin{eqnarray}\label{e15}
&&\sum\limits_{r=1}^{m_{n}}\sum\limits_{i=1}^{j+l_{r}}E(\left|\xi_{i,r}\right|^{3}e^{t\left|\xi_{i,r}\right|})\nonumber\\
&\leq &\Delta_{j,1}^{-3/2} \sum\limits_{r=1}^{m_{n}}d_{j,r}^{2}\Big(\frac{6e^{-tD_{j,r}}}{(1-D_{j,r}t)^{4}}+\frac{e^{-tD_{j,r}}}{1-D_{j,r}t}
+\frac{6e^{tD_{j,r}}}{(1+D_{j,r}t)^{4}}+\frac{e^{tD_{j,r}}}{1+D_{j,r}t}\Big).
\end{eqnarray}
The above estimate is valid if $tD_{j,r}<1$ for all $n-j_n+1\le j\le n$ and $1\le r\le m_n$.

When $n-j_n+1\le j\le n$ and $1\le r\le m_n$, we have $j+l_{r}>n-j_{n}\geq n/2$, $\Delta_{j,1}\ge 1/(j+l_1)=1/j\ge 1/n$, and
$d_{j, r}=\frac{1}{j+l_{r}}\leq 2/n$. Therefore,
\[
D_{j,r}=\frac{d_{j,r}}{\Delta_{j,1}^{1/2}}\le \frac{2}{n^{1/2}},
\]
which implies
\[
tD_{j,r}\le 2t_nn^{-1/2}=O(n^{-5/14})\to 0
\]
uniformly over $0\le t\le t_n=O(n^{1/7})$,
 $n-j_n+1\le j\le n$ and $1\le r\le m_n$ as $n\to\infty$.  Hence, it follows from \eqref{e15} and Lemma~\ref{l6} that for some constant $C>0$
\begin{equation}\label{e16}
\gamma:=\sum\limits_{r=1}^{m_{n}}\sum\limits_{i=1}^{j+l_{r}}E(\left|\xi_{i,r}\right|^{3}e^{t\left|\xi_{i,r}\right|})\leq \frac{C\sum\limits_{r=1}^{m_{n}}d_{j,r}^{2}}{\Delta_{j,1}^{3/2}}=\frac{C\Delta_{j,2}}{\Delta_{j,1}^{\frac{3}{2}}}\le \frac{C}{j^{1/2}}\le \frac{2C}{n^{1/2}}
\end{equation}
uniformly over $n-j_n+1\le j\le n$ as $n\to\infty$.

By Lemma~\ref{l9}, $\frac{P(W_{j}\geq t)}{1-\Phi(t)}=1+O(1)(1+t^{3})\gamma e^{4t^{3}\gamma}=1+O(n^{-1/14})$ uniformly over
$0\leq t\leq t_n$ and  $n-j_{n}+1\leq j\leq n$ as $n\to\infty$. \hfill$\blacksquare$

\subsection{Proofs of Theorems~\ref{t1} and \ref{t2}}\label{proof2}

\noindent {\it Proof of Theorem~\ref{t1}}.  Define
\begin{equation}\label{jn}
j_{n}=\mbox{the integer part of } \Delta_n^{-1/2}\cdot n^{1/7}+1.
\end{equation}

The proof of the theorem will be divided into three steps.

{\it Step 1.} We will prove that
\begin{equation}\label{e18}\lim\limits_{n\rightarrow\infty}\sum\limits_{j=1}^{n-j_{n}}P(\ln G_{j}>\sum\limits_{r=1}^{m_{n}}\psi(n+l_{r})+\Delta_n^{1/2}(a_{n}+b_{n}y))=0,\ y\in\mathbb{R}.
\end{equation}

Since $\Delta_n\ge 1/n$, we have from \eqref{jn} that
\[
\frac{j_{n}}{n}\le \frac{n^{1/7}}{n\Delta_n^{1/2}}+\frac{1}{n}\le \frac{2}{n^{5/14}}\to 0
\]
and
\[
j_{n}\big(\frac{\Delta_n}{\ln n}\big)^{1/2}\ge \frac{n^{1/7}}{\Delta_n^{1/2}}\frac{\Delta_n^{1/2}}{(\ln n)^{1/2}}=\frac{n^{1/7}}{(\ln n)^{1/2}}\to\infty,
\]
as $n\to\infty$, that is, the conditions in Lemma~\ref{l8} are satisfied.  Therefore, (\ref{e18}) holds in case $\alpha\in (0, \infty]$. In case $\alpha=0$, $a_{n}+b_{n}y>0$ for all large $n$, by Lemma~\ref{l8}, we have
\begin{eqnarray*}
& &\lim\limits_{n\rightarrow\infty}\sum\limits_{j=1}^{n-j_{n}}P\big(\ln G_{j}>\sum\limits_{r=1}^{m_{n}}\psi(n+l_{r})+\Delta_n^{1/2}(a_{n}+b_{n}y)\big)\\
&\leq& \lim\limits_{n\rightarrow\infty}\sum\limits_{j=1}^{n-j_{n}}P\big(\ln G_{j}>\sum\limits_{r=1}^{m_{n}}\psi(n+l_{r})\big)\\
&=&0.
\end{eqnarray*}
Note that \eqref{e18} implies
\[
\lim_{n\to\infty}P\Big(\frac{\max\limits_{1\leq j\leq n-j_{n}}\ln G_{j} -\sum\limits_{r=1}^{m_{n}}\psi(n_{r})}{\Delta_n^{1/2}b_{n}}-\frac{a_{n}}{b_{n}}> y\Big)=0, ~~~~y\in \mathbb{R}
\]
or equivalently
\begin{equation}\label{e19}
\lim_{n\to\infty}P\Big(\frac{\max\limits_{1\leq j\leq n-j_{n}}\ln G_{j} -\sum\limits_{r=1}^{m_{n}}\psi(n_{r})}{\Delta_n^{1/2}b_{n}}-\frac{a_{n}}{b_{n}}\le y\Big)=1, ~~~~y\in \mathbb{R}.
\end{equation}

\medskip

{\it Step 2.} We claim that
\begin{equation}\label{e20}
\frac{M_{n}(j_{n})}{\Delta_n^{1/2}b_{n}} ~~~\mbox{converges in probability to zero}.
\end{equation}
To prove this, it suffices to show that $M_{n}(j_{n})=O_{p}(\Delta_n^{1/2}(\ln n)^{-1})$ since $b_{n}\geq(\ln n)^{-1/2}$ for large $n$.

When $\alpha\in\left( 0,\infty\right]$, $\Delta_n^{-1/2}$ is bounded, and $j_{n}=O(n^{1/7})$.  By Lemma~\ref{l7}, we have
\[
M_{n}(j_{n})=O_{p}(j_{n}(\frac{\Delta_n}{n})^{1/2})=O_{p}(\Delta_n^{1/2}n^{-5/15})=O_{p}(\Delta_n^{1/2}(\ln n)^{-1}).
\]

When $\alpha=0$,  by Lemma~\ref{l7}, we can obtain that
\begin{eqnarray*}
M_{n}(j_{n})&=&O_{p}(\min\left\lbrace j_{n}(\frac{\Delta_n}{n})^{1/2},\Delta_n\ln n\right\rbrace )\\
&=&\Delta_n^{1/2}O_{p}(\min\left\lbrace \Delta_n^{-1/2}n^{-5/14},\Delta_n^{1/2}\ln n\right\rbrace )\\
&=&\Delta_n^{1/2}\cdot O_{p}(n^{-1/8})\\
&=&O_{p}(\Delta_n^{1/2}(\ln n)^{-1/2})
\end{eqnarray*}
since $\Delta_n^{-1/2}n^{-5/14}\leq n^{-1/8}$ if $\Delta_n^{-1/2}\leq n^{1/7}$ and $\Delta_n^{1/2}\ln n\leq n^{-1/8}$
if $\Delta_n^{-1/2}> n^{1/7}$.  This proves \eqref{e20}.

\medskip

{\it Step 3.} Set
\[
T_{n}(j_{n})=\max\limits_{n-j_{n}+1\leq j\leq n}\Big\{\sum\limits_{r=1}^{m_{n}}\frac{s_{j,r}-(j+l_{r})}{j+l_{r}}+\sum\limits_{r=1}^{m_{n}}{\psi(j+l_{r})}\Big\}.
\]
We will show that for every $y\in \mathbb{R}$
\begin{equation}\label{e21}
P\big(T_n(j_n)\le\sum\limits_{r=1}^{m_{n}}\psi(n_{r})+\Delta_n^{1/2}(a_{n}+b_{n}y)\big)\rightarrow\Phi_{\alpha}(y).
\end{equation}

In fact,
\begin{eqnarray}\label{product}
&&P\big(T_{n}(j_{n})\leq\sum\limits_{r=1}^{m_{n}}\psi(n+l_{r})+\Delta_n^{1/2}(a_{n}+b_{n}y)\big)\nonumber\\
&=&\prod_{j=n-j_{n}+1}^{n}P\big(W_{j}\leq\frac{1}{\Delta_{j,1}^{1/2}}
(\sum\limits_{r=1}^{m_{n}}(\psi(n+l_{r})-\psi(j+l_{r}))+\Delta_n^{1/2}(a_{n}+b_{n}y))\big)\nonumber\\
&=&\prod_{i=1}^{j_{n}}P\Big(W_{n-i+1}\leq \frac{\sum\limits_{r=1}^{m_{n}}(\psi(n_{r})-\psi(n_{r}-i+1))+\Delta_n^{1/2}(a_{n}+b_{n}y)}{\big(\sum\limits_{r=1}^{m_{n}}\frac{1}{n_{r}-i+1}\big)^{1/2}}
\Big)\nonumber\\
&=&\prod\limits_{i=1}^{j_{n}}(1-a_{ni}),
\end{eqnarray}
where $a_{ni}=P(W_{n-i+1}\geq t_{n,i})$ and
$$
t_{n,i}=\big(\sum\limits_{r=1}^{m_{n}}\frac{1}{n_{r}-i+1}\big)^{-1/2}
\big(\sum\limits_{r=1}^{m_{n}}(\psi(n_{r})-\psi(n_{r}-i+1))+\Delta_n^{1/2}(a_{n}+b_{n}y)\big).
$$
It follows from \eqref{e9} and Taylor's expansion that
\begin{eqnarray*}
&&
\big(\sum_{r=1}^{m_{n}}\frac{1}{n_{r}-i+1}\big)^{-1/2}\sum_{r=1}^{m_{n}}\big(\psi(n_{r})-\psi(n_{r}-i+1)\big)\\
&=&
\big(\sum\limits_{r=1}^{m_{n}}\frac{1}{n_{r}}\cdot\frac{1}{1-\frac{i-1}{n_{r}}}\big)^{-1/2}\sum\limits_{r=1}^{m_{n}}\frac{i-1}{n_{r}}(1+O(\frac{i}{n_{r}}))\\
&=&
(i-1)\big(\sum\limits_{r=1}^{m_{n}}\frac{1}{n_{r}}(1+O(\frac{i-1}{n_{r}}))\big)^{-1/2}\sum\limits_{r=1}^{m_{n}}\frac{1}{n_{r}}(1+O(\frac{i}{n_{r}}))\\
&=&
(i-1)\big(1+O(\frac{j_n}{n})\big)(\sum\limits_{r=1}^{m_{n}}\frac{1}{n_{r}})^{1/2}\\
&=&
(i-1)\big(1+O(n^{-5/14})\big)\Delta_n^{1/2}
\end{eqnarray*}
and
 \begin{eqnarray*}
&&\Big( \big(\sum\limits_{r=1}^{m_{n}}\frac{1}{n_{r}-i+1}\big)^{-1/2}\Delta_n^{1/2}-1\Big)(a_{n}+b_{n}y)\\
&=&
\Big( \big(\sum\limits_{r=1}^{m_{n}}\frac{1}{n_{r}}(1+O(\frac{i-1}{n_{r}}))\big)^{-1/2}\cdot \Delta_n^{1/2}-1\Big) (a_{n}+b_{n}y)\\
&=&\Big(\big(\sum\limits_{r=1}^{m_{n}}\frac{1}{n_{r}}+O(\sum^{m_n}_{r=1}\frac{1}{n_r^2}))\big)^{-1/2}\cdot \Delta_n^{1/2}-1\Big) (a_{n}+b_{n}y)\\
&=&\Big(\big(\Delta_n+O(\Delta_{n,2})(i-1)\big)^{-1/2}\Delta_n^{1/2}-1\Big) (a_{n}+b_{n}y)\\
&=&\Big(\big(1+O(\frac{\Delta_{n,2}}{\Delta_n})(i-1)\big)^{-1/2}-1\Big) (a_{n}+b_{n}y)\\
&=&
O\Big(\frac{(i-1)\Delta_{n,2}}{\Delta_n}(\ln n)^{1/2}\Big)\\
&=&\Delta_n^{1/2}(i-1)\cdot O(\frac{\Delta_{n,2}}{\Delta_n^{1.5}}(\ln n)^{1/2})\\
&=&\Delta_n^{1/2}(i-1)\cdot O(\frac{(\ln n)^{1/2}}{n^{1/2}}).
\end{eqnarray*}
In the above estimation we have used the facts (a): $\max_{1\le i\le j_n}(i-1)\Delta_{n,2}/\Delta_n\le j_n/n\to 0$ from Lemma~\ref{l6};
(b): $a_n+b_ny=O((\ln n)^{1/2})$; and (c): $\Delta_{n,2}/\Delta_n^{1.5}\le n^{-1/2}$ from Lemma~\ref{l6}.
Therefore, we conclude that
\begin{equation}\label{e22}
t_{n,i}=(i-1)\big(1+O(n^{-5/14})\big)\Delta_n^{1/2}+a_{n}+b_{n}y
\end{equation}
holds uniformly over $1\le i\le j_n$ as  $n\to\infty$.

\medskip

{\it Case 1.} If $\alpha=0$, then $\Delta_n\to 0$ and
\[
a_{n}=a(\Delta_n^{-1})\sim(\ln(\Delta_n^{-1}))^{1/2} ~~\mbox{ and }~ b_{n}=b(\Delta_n^{-1})\sim(\ln(\Delta_n^{-1}))^{-1/2},
\]
we have
$$\min\limits_{1\leq i\leq j_{n}}t_{n,i}\rightarrow\infty~~\mbox{ and }~~\max\limits_{1\leq i\leq j_{n}}t_{n,i}=O(\Delta_n^{1/2}j_{n}+(\ln n)^{1/2})=O(n^{\frac{1}{7}}).
$$
It follows from Lemma~\ref{l10} that
\begin{equation}\label{ani}
a_{ni}=(1+o(1))(1-\Phi(t_{n,i}))
\end{equation}
uniformly over $1\leq i\leq j_{n}$.

Now define $c_{n,i}$ such that $t_{n,i}=(i-1)c_{n,i}+a_n+b_ny$ with $c_{n,1}=0$ and apply Lemma~\ref{l4} with  $x_{n}=\Delta_n^{-1}$ by noting that $c_{n,i}=\big(1+O(n^{-5/14})\big)\cdot\Delta_n^{1/2}$ from \eqref{e22}. Then we get
$$
\sum\limits_{i=1}^{j_{n}}a_{ni}=(1+o(1))\sum\limits_{i=1}^{j_{n}}(1-\Phi(t_{n,i}))\rightarrow e^{-y}.
$$
It is obvious from \eqref{ani} that $\max\limits_{1\leq i\leq j_{n}}a_{ni}\rightarrow0$. So we have from Lemma~\ref{prod2sum} that
$\prod^{j_n}_{i=}(1-a_{ni})\to \exp(-e^{-y})=\Phi_0(y)$ as $n\to\infty$, which together with \eqref{product} yields \eqref{e21} with $\alpha=0$,

\medskip

{\it Case 2.} If $\alpha\in(0,\infty)$, then $j_{n}\sim\alpha^{-1/2}n^{1/7}$. By definition, $a_{n}=0$ and $b_{n}=1$, and \eqref{e22} means
$$
t_{n,i}=(1+o(1))\alpha^{1/2}(i-1)+y
$$
holds uniformly over $1\leq j\leq j_{n}$ as $n\to\infty$.

Let $j_{0}>1$ be an integer such that $\min\limits_{j_{0}\leq i\leq j_{n}}t_{n,i}>0$. Since $\max\limits_{1\leq i\leq j_{n}}\left| t_{n,i}\right|=O(n^{1/7})$, we have from Lemma~\ref{l10}
\begin{equation}\label{e23}a_{ni}=(1+o(1))(1-\Phi(t_{n,i}))
\end{equation}
uniformly over $j_{0}\leq i\leq j_{n}$. By using the standard central limit theorem, we know this also holds for each $i=1,2,\cdots,j_{0}-1$. Therefore, for each $i\geq1$,
\begin{equation}\label{e24}
\lim\limits_{n\to\infty}a_{ni}=1-\Phi(\alpha^{1/2}(i-1)+y)
\end{equation}
and
\begin{equation}
\label{e25}\sum\limits_{i\geq1}(1-\Phi(\alpha^{1/2}(i-1)+y))<\infty\end{equation} by the fact $1-\Phi(x)\sim\frac{1}{\sqrt{2\pi}x}e^{-x^2/2}$ as $x\to+\infty$.

Define $a_{n,i}=0$ for $i>j_{n}$. By the fact that $t_{n,i}\geq \frac12\alpha^{1/2}(i-1)+y\geq y$ for $1\leq i\leq j_{n}$ for all large $n$, we have $\sup_{n\geq n_{0},1\leq i\leq j_{n}}a_{ni}<1$ for some integer $n_{0}$. And since $a_{ni}\leq2(1-\Phi(\frac12\alpha^{1/2}(i-1)+y))$ for all $1\leq i\leq j_{n}$ as $n$ is sufficiently large and $\sum\limits_{i\geq1}2(1-\Phi(\alpha^{1/2}(i-1)+y))<\infty$, we obtain that $\lim\limits_{n\to\infty}\sum\limits_{i=1}^{j_{n}}a_{ni}=\sum\limits_{i=1}^{\infty}(1-\Phi(\alpha^{1/2}(i-1)+y))$. So it follows from Lemma~\ref{l3} that \[
\lim\limits_{n\to\infty}\prod\limits_{i=1}^{j_{n}}(1-a_{ni})=\prod\limits_{i=1}^{\infty}\Phi(y+\alpha^{1/2}(i-1))=\Phi_{\alpha}(y),
\]
which together with \eqref{product} yields \eqref{e21} with $\alpha\in (0, \infty)$.

\medskip

{\it Case 3.} If $\alpha=\infty$, then by the fact $0\leq \Delta_n^{1/2}(j_{n}-1)\leq n^{1/7}$, we have $t_{n,i}=O(n^{1/7})$. In particular, we have $t_{n,1}=y$ and for all large $n$, $t_{n,i}>0$ if $2\leq i\leq j_{n}$ and $j_{n}\geq2$. So we obtain  from Lemma~\ref{l10} that
$$a_{ni}=(1+o(1))(1-\Phi(t_{n,i}))$$
uniformly over $1\leq i\leq j_{n}$. Note that $t_{n,i}\geq\frac{i}{3}\Delta_n^{1/2}$ if $2\leq i\leq j_{n}$ and $j_{n}\geq2$. For large $n$ we have
\[
I(j_{n}\geq2)\sum\limits_{i=2}^{j_{n}}t_{n,i}\leq2\sum\limits_{i=2}^{\infty}(1-\Phi(\frac{i}{3}\Delta_n^{1/2}))
\leq\sum\limits_{i=2}^{\infty}\exp(-\frac{i^{2}}{18}\Delta_n)
\leq3\sqrt{2}\pi\Delta_n^{-1/2}\to 0
\]
since $\exp(-\frac{i^{2}}{18}\sum\limits_{r=1}^{m_{n}}\frac{1}{n_{r}})\leq\int_{i-1}^{i}\exp(-\frac{x^{2}}{18}\sum\limits_{r=1}^{m_{n}}\frac{1}{n_{r}})dx$ for $i\geq2$. It is also obvious that $I(j_{n}\geq2)\max\limits_{2\leq i\leq j_{n}}a_{ni}\to 0$, so $I(j_{n}\geq2)(1-\prod\limits_{i=2}^{j_{n}}(1-a_{ni}))\to0$ as $n\to\infty$, which coupled with \eqref{product} implies
\begin{eqnarray*}
&&P(T_{n}(j_{n})\leq\sum\limits_{r=1}^{m_{n}}\psi(n_{r})+\Delta_n^{1/2}(a_{n}+b_{n}y))\\
&=&\prod\limits_{i=1}^{j_{n}}(1-a_{ni})\\
&=&(1-a_{n1})\Big(1-I(j_{n}\geq2)\big(1-\prod\limits_{i=2}^{j_{n}}(1-a_{ni})\big)\Big)\\
&\to&\Phi(y)=\Phi_{\infty}(y),
\end{eqnarray*}
i.e. \eqref{e21} holds with $\alpha=\infty$.

\medskip

Now we are ready to conclude the proof. We first have from \eqref{e21} that
\[
\frac{T_n(j_n)-\sum\limits_{r=1}^{m_{n}}\psi(n_{r})}{\Delta_n^{1/2}b_n}-\frac{a_n}{b_n}\xrightarrow{d}\Phi_{\alpha}.
\]
By Lemma~\ref{l5} and \eqref{e20}, we get
\[
\frac{\max\limits_{n-j_{n}+1\leq j\leq n}\ln G_{j} -\sum\limits_{r=1}^{m_{n}}\psi(n_{r})}{\Delta_n^{1/2}b_{n}}-\frac{a_{n}}{b_{n}}\xrightarrow{d}\Phi_{\alpha},
\]
or equivalently
\[
\lim_{n\to\infty}P\Big(\frac{\max\limits_{n-j_{n}+1\leq j\leq n}\ln G_{j} -\sum\limits_{r=1}^{m_{n}}\psi(n_{r})}{\Delta_n^{1/2}b_{n}}-\frac{a_{n}}{b_{n}}\le y\Big)=\Phi_{\alpha}(y),~~~y\in\mathbb{R},
\]
which together with \eqref{e19} and the independence of $\max\limits_{1\leq j\leq n-j_n}\ln G_{j}$ and $\max\limits_{n-j_n+1\le j\le n}\ln G_{j}$ yields that
\begin{eqnarray*}
&&P\Big(\frac{\max\limits_{1\leq j\leq n}\ln G_{j}-\sum\limits_{r=1}^{m_{n}}\psi(n_{r})}{\Delta_n^{1/2}b_{n}}-\frac{a_{n}}{b_{n}}\le y\Big)\\
&=&P\Big(\frac{\max\limits_{1\leq j\leq n-j_n}\ln G_{j} -\sum\limits_{r=1}^{m_{n}}\psi(n_{r})}{\Delta_n^{1/2}b_{n}}-\frac{a_{n}}{b_{n}}\le y\Big)\\
&&~~~\times
P\Big(\frac{\max\limits_{n-j_{n}+1\leq j\leq n}\ln G_{j}-\sum\limits_{r=1}^{m_{n}}\psi(n_{r})}{\Delta_n^{1/2}b_{n}}-\frac{a_{n}}{b_{n}}\le y\Big)\\
&\to&\Phi_{\alpha}(y)
\end{eqnarray*}
for every $y\in \mathbb{R}$. Since $G_j=\prod^{m_n}_{r=1}s_{j,r}$,  $\max\limits_{1\leq j\leq n}\ln |\bz_{j}|$ and $\frac12\max\limits_{1\leq j\leq n}\ln G_{j}$
have the same distribution from Lemma~\ref{p1}. Hence we conclude that
 \[
\lim\limits_{n\rightarrow\infty}
P\Big(\frac{\max\limits_{1\leq j\leq n}\ln |\bz_{j} | -\sum\limits_{r=1}^{m_{n}}\psi(n_{r})/2}{\Delta_n^{1/2}/2}\leq a_{n}+b_{n}y\Big)=\Phi_{\alpha}(y),
\]
proving \eqref{e17}. This completes the proof of Theorem~\ref{t1}.\hfill$\blacksquare$\\

\noindent{\it Proof of Theorem~\ref{t2}}.  Define for $\alpha\in[0,\infty)$,
$$V_{n}=\frac{\max\limits_{1\leq j\leq n}\ln\left| \bz_{j}\right|-\sum\limits_{r=1}^{m_{n}}\psi(n_{r})/2}{\Delta_n^{1/2}b_{n}/2}-\frac{a_{n}}{b_{n}}.$$
Then $V_{n}$ converges in distribution to $\Theta_{\alpha}$, where $\Theta_{\alpha}$ is a random variable with the cdf $\Phi_{\alpha}(y)$. And it can be easily verified that
\begin{align}\begin{split}\label{e28}
\max\limits_{1\leq j\leq n}\left| \bz_{j}\right| &=\exp\left\lbrace \frac12\sum\limits_{r=1}^{m_{n}}\psi(n_{r})+\frac12\Delta_n^{1/2}(a_{n}+b_{n}V_{n})\right\rbrace\\&=\exp\left\lbrace \frac12\sum\limits_{r=1}^{m_{n}}\psi(n_{r})+\frac12\Delta_n^{1/2}a_{n}\right\rbrace \cdot\exp\left\lbrace \frac12\Delta_n^{1/2}b_{n}V_{n}\right\rbrace. \end{split}\end{align}

\medskip

(a). If $\alpha=0$, then we have $\Delta_n\to0$, $a_{n}=a(x_{n})\sim(\ln\Delta_n^{-1})^{1/2}\to\infty$, $b_{n}=b(\Delta_n^{-1})\sim(\ln\Delta_n^{-1})^{-1/2}\to 0$, and $\Delta_n^{1/2}a_{n}\sim\Delta_n^{1/2}b_{n}^{-1}$ as $n\to\infty$.
Thus, we get from \eqref{e9} and Taylor's expansion that
\begin{eqnarray*}
\max\limits_{1\leq j\leq n}\left| \bz_{j}\right|
&=&\exp\Big\{ \frac12\sum\limits_{r=1}^{m_{n}}\ln n_{r}+O(\Delta_n)+\frac12\Delta_n^{1/2}a_{n}\Big\}\cdot (1+\frac12\Delta_n^{1/2}b_{n}V_{n}+O_{p}(b_{n}^{2}\Delta_n))\\
&=&\big(\prod\limits_{r=1}^{m_{n}}n_{r}\big)^{1/2}(1+\frac12\Delta_n^{1/2}a_{n}+O(\Delta_n))(1+\frac12\Delta_n^{1/2}b_{n}V_{n}+O_{p}(\Delta_n))\\
&=&\big(\prod\limits_{r=1}^{m_{n}}n_{r}\big)^{1/2}(1+\frac12\Delta_n^{1/2}a_{n}+\frac12\Delta_n^{1/2}b_{n}V_{n}+O_{p}(\Delta_n a_{n}^{2})),
\end{eqnarray*}
which implies that
$$\frac{1}{\Delta_n^{1/2}b_{n}/2}\Big(\frac{\max\limits_{1\leq j\leq n}|\bz_{j}| }{(\prod\limits_{r=1}^{m_{n}}n_{r})^{1/2}}-1\Big)-\frac{a_{n}}{b_{n}}=V_{n}+O_{p}(\Delta_n^{1/2}(\ln\Delta_n^{-1})^{3/2})
$$ converges in distribution to $\Lambda$.\\

(b). If $\alpha\in(0,\infty)$, then $a_{n}=0$ and $b_{n}=1$ in this case. Therefore, we have
\begin{eqnarray*}
\max\limits_{1\leq j\leq n} | \bz_{j}| &=&\exp\Big\{ \frac12\sum\limits_{r=1}^{m_{n}}\psi(n_{r})+\frac12\Delta_n^{1/2}V_{n}\Big\}\\
&=&\exp\Big\{\frac12\sum\limits_{r=1}^{m_{n}}\psi(n_{r})\Big\} \cdot\exp\Big\{ \frac12\Delta_n^{1/2}V_{n}\Big\}.
\end{eqnarray*}
Using \eqref{e9}, we have  $\sum\limits^{m_n}_{r=1}\psi(n_{r})=\sum^{m_n}_{r=1}\ln n_{r}-\frac{1}{2}\Delta_n+o(\Delta_n)$, and then we obtain
$$
\frac{\max\limits_{1\leq j\leq n}| \bz_{j}|}{\big(\prod\limits_{r=1}^{m_{n}}n_{r}\big)^{1/2}}
=\exp\big(-\frac{1}{4}\alpha+o(1)\big)\cdot\exp\big((\frac12\alpha^{1/2}+o(1))V_{n}\big),
$$
which converges in distribution to $\Phi_{\alpha}(\frac12\alpha^{1/2}+2\alpha^{-1/2}\ln y)$, $y>0$, the cumulative
distribution of $e^{-\alpha/4}\exp\big(\frac12\alpha^{1/2}\Theta_{\alpha}\big)$. This completes the proof.
\hfill$\blacksquare$

\vspace{10pt}

\noindent\textbf{Acknowledgements}. The authors would like to thank an anonymous referee for his/her careful reading of the manuscript.
The research of Yongcheng Qi was supported in part by NSF Grant DMS-1916014.


\begin{thebibliography}{AA} 


\bibitem{Abramowitz1972} Abramowitz, M. and Stegun, I. A. (1972). Handbook of Mathematical Functions. {\it Dover, New York}.


\bibitem{adhi2016}
Adhikari, K., Reddy,  N. K., Reddy, T. R. and Saha, K. (2016). Determinantal point processes in the plane
from products of random matrices. {\it Ann. Inst. H. Poincare Probab. Statist.} 52 (1), 16-46.


\bibitem{akemann2011oxford}
Akemann, G., Baik, J. and Di Francesco, P. (2011). The Oxford handbook of random matrix theory. {\it Oxford University Press}.

\bibitem{akemann2013products}
Akemann, G., Ipsen, J.R. and Kieburg, M. (2013). Products of rectangular random matrices: singular values and progressive scattering. {\it Physical Review} E, 88 (5), 052118.



\bibitem{Bai}
Bai, Z. D. (1999). Methodologies in spectral analysis of large dimensional random matrices,
a review. {\it Statistica Sinica} 9, 9611-9677.


\bibitem{Baik1999}
Baik, J, Deift, P. and Johansson, K. (1999). On the distribution of the length of the longest increasing
subsequence of random permutations. {\it Jour. Amer. Math. Soc.} 12(4), 1119-1178.


\bibitem{Been1997}
Beenakker, C. W. J. (1997). Random-matrix theory of quantum transport.
{\it Rev. Mod. Phys.} 69, 731-809.


\bibitem{benet2014spectral}
Benet, L. (2014). Spectral domain of large nonsymmetric correlated Wishart matrices. {\it Physical Review} E, 90(4), 042109.


\bibitem{Bor}
Bordenave, C. (2011). On the spectrum of sum and product of non-Hemitian random matrices. {\it Elect. Comm. in Probab.} 16, 104-113.


\bibitem{bouchaud2007large}
Bouchaud, J.P., Laloux, L., Miceli, M.A. and Potters, M. (2007). Large dimension forecasting models and random singular value spectra. {\it The European Physical Journal B} 55(2),  201-207.

\bibitem{bouchaud2009financial}
Bouchaud, J.P. and Potters, M. (2009). Financial applications of random matrix theory: a short review. {\it arXiv preprint arXiv:0910.1205}.


\bibitem{Burda}
Burda, Z. (2013).  Free products of large random matrices - a short review of recent developments. {\it J. Phys. Conf. Ser. 473, 012002}. Also available at {\tt http://arxiv.org/pdf/1309.2568v2.pdf}.

\bibitem{BJW}
Burda, Z., Janik, R. A. and Waclaw, B. (2010).  Spectrum of the product of independent random Gaussian matrices.  {\it Phys. Rev.} E 81, 041132.





\bibitem{ChangLiQi2018}
Chang, S., Li, D. and Qi, Y. (2018).  Limiting distributions of spectral radii for product of matrices from the spherical ensemble.
 {\it Journal of Mathematical Analysis and Applications} 461, 1165-1176.



\bibitem{ChangQi2017}
Chang, S. and Qi, Y. (2017). Empirical distribution of scaled eigenvalues for product of matrices from the spherical ensemble.
{\it Statistics and Probability Letters} 128, 8-13.


\bibitem{chen2013}
Chen, L.H., Fang, X. and Shao, Q.M. (2013). From Stein identities to moderate deviations. {\it The Annals of Probability} 41(1),  262-293.

\bibitem{chow2012}
Chow, Y.S. and Teicher, H. (2012). Probability Theory: Independence, Interchangeability, Martingales. {\it Springer Science and Business Media}.

\bibitem{CD}
Couillet, R.  and Debbah, M. (2011). {\it Random matrix methods for wireless communications}. Cambridge
Univ Press.


\bibitem{crisanti1994products}
Crisanti, A., Paladin, G. and Vulpiani, A. (2012). Products of Random Matrices: in Statistical Physics (Vol. 104). {\it Springer Science and Business Media}.




\bibitem{Goetz}
G\"{o}tze, F. and Tikhomirov, T. (2010). On the asymptotic spectrum of products of independent random
matrices. {\tt http://arxiv.org/pdf/1012.2710v3.pdf}.

\bibitem{GuiQi2018}
Gui, W. and Qi, Y. (2018). Spectral radii of truncated circular unitary matrices. {\it Journal of Mathematical Analysis and Applications} 458(1),  536-554.

\bibitem{haake1991quantum}
Haake, F. (2013). Quantum signatures of chaos (Vol. 54). {\it Springer Science and Business Media}.



\bibitem{Ipsen2015}
Ipsen, J.R. (2015). Products of independent Gaussian random matrices. {\it Doctoral Dissertation}, Bielefeld University.


\bibitem{ipsen2014weak}
Ipsen, J.R. and Kieburg, M. (2014). Weak commutation relations and eigenvalue statistics for products of rectangular random matrices. {\it Physical Review} E, 89(3), 032106.



\bibitem{Jiang09}
Jiang, T.  (2009). Approximation of Haar distributed matrices and limiting distributions of eigenvalues of Jacobi ensembles.  {\it Probability Theory and Related Fields} 144 (1), 221-246.



\bibitem{JiangQi2017}
Jiang, T. and Qi, Y. (2017). Spectral radii of large non-Hermitian random matrices. {\it Journal of Theoretical Probability} 30(1),  326-364.

\bibitem{JiangQi2018}
Jiang, T. and Qi, Y. (2019). Empirical distributions of eigenvalues of product ensembles. {\it Journal of Theoretical Probability} 32(1), 353-394.


\bibitem{Johansson07}
Johansson, K. (2007). From Gumbel to Tracy-Widom. {\it Probab. Theory Relat. Fields} 138, 75-112.

\bibitem{John2001}
Johnstone, I. (2001). On the distribution of the largest eigenvalue in principal components analysis. {\it Ann. Stat.} 29, 295-327.

\bibitem{John2008}
Johnstone, I. (2008). Multivariate analysis and Jacobi ensembles: Largest eigenvalue, Tracy-Widom limits and rates of convergence. {\it Ann. Stat.}, 36 (6), 2638-2716.


\bibitem{La2018}
Lacroix-A-Chez-Toine, B.,  Grabsch, A.,  Majumdar,S. N. and  Schehr G. (2018). Extreme of 2d Coulomb gas: universal intermediate deviation regime.
{\it Journal of Statistical Mechanics: Theory and Experiment} 013203. DOI: 10.1088/1742-5468/aa9bb2.

\bibitem{livan2012asymmetric}
Livan, G. and Rebecchi, L. (2012). Asymmetric correlation matrices: an analysis of financial data. {\it The European Physical Journal} B, 85(6), 213.



\bibitem{MS2005}
Mezzadri, F. and Snaith, N. C. (2005). {\it Recent perspectives in random matrix theory and number theory}.
Cambridge Univ Press.


\bibitem{muller2002asymptotic}
Muller, R.R. (2002). On the asymptotic eigenvalue distribution of concatenated vector-valued fading channels. {\it IEEE Transactions on Information Theory} 48(7), 2086-2091.


\bibitem{Rourke}
O'Rourke, S. and Soshnikov, A.  (2011). Products of independent non-Hermitian random matrices. {\it Electrical Journal of Probability} 16(81), 2219-2245.

\bibitem{Rourke14}
O'Rourke, S., Renfrew, D., Soshnikov, A. and Vu, V.  (2014). Products of independent elliptic random matrices. Available at {\tt http://arxiv.org/pdf/1403.6080v2.pdf}.


\bibitem{rider2003limit}
Rider, B. C. (2003). A limit theorem at the edge of a non-Hermitian random matrix ensemble. {\it Journal of Physics} A: Mathematical and General, 36(12), 3401-3409.

\bibitem{rider2004order}
Rider, B. C. (2004). Order statistics and Ginibre's ensembles. {\it Journal of statistical physics} 114(3-4),  1139-1148.

\bibitem{rider2014extremal}
Rider, B. C. and Sinclair, C.D. (2014). Extremal laws for the real Ginibre ensemble. {\it The Annals of Applied Probability} 24(4),  1621-1651.


\bibitem{Tracy94}
Tracy, C. A. and Widom, H. (1994). Level-spacing distributions and Airy kernal. {\it Comm.
Math. Physics} 159, 151-174.


\bibitem{Tracy96}
Tracy, C. A. and Widom, H. (1996). On the orthogonal and symplectic matrix ensembles.
{\it Comm. Math. Physics} 177, 727-754.

\bibitem{TW02}
Tracy, C. A. and Widom, H. (2002).  Distribution functions for largest eigenvalues and their applications. {\it Proceedings of the
ICM, Beijing} 1, 587-596.


\bibitem{tulino2004random}
Tulino, A.M. and Verd\'u, S. (2004). Random matrix theory and wireless communications. {\it Foundations and Trends® in Communications and Information Theory} 1(1),  1-182.

\bibitem{vinayak2013spectral}
Vinayak (2013). Spectral density of a Wishart model for nonsymmetric correlation matrices. {\it Physical Review} E, 88(4), 042130.

\bibitem{Wigner}
Wigner, E. P. (1955). Characteristic vectors of bordered matrices with infinite dimensions. {\it Ann. Math.} 62, 548-564.

\bibitem{Wishart}
Wishart, J. (1928). The generalized product moment distribution in samples
from a normal multivariate population. {\it Biometrika} 20, 35-52.

\bibitem{Zeng2016}
Zeng, X. (2016). Eigenvalues distribution for products of independent spherical ensembles. {\it J. Phys. A: Math. Theor.} 49, 235201.


\bibitem{Zeng2017}
Zeng, X. (2017). Limiting empirical distribution for eigenvalues of products of random rectangular matrices. {\it Statistics and Probability Letters} 126, 33-40.


\end{thebibliography}

\baselineskip 12pt
\def\ref{\par\noindent\hangindent 25pt}

\end{document}